\documentclass{aims}

\usepackage{txfonts}
\def\typeofarticle{Research article}
\def\currentvolume{0}
\def\currentissue{0}
\def\currentyear{2026}

\def\ppages{0--0}
\def\DOI{10.3934/math.2026xxx}
\def\Received{XX XXX 2026}
\def\Revised{XX XXX 2026}
\def\Accepted{XX XXX 2026}
\def\Published{XXX 2026}

\usepackage{amsmath, amssymb, amsthm}
\usepackage{mathrsfs}
\usepackage{graphicx}

\numberwithin{equation}{section}

\newcommand{\R}{\mathbb{R}}
\newcommand{\N}{\mathbb{N}}
\newcommand{\norm}[1]{\lVert#1\rVert}
\newcommand{\diam}{\operatorname{diam}}

\newtheorem{theorem}{Theorem}[section]
\newtheorem{lemma}[theorem]{Lemma}
\newtheorem{definition}[theorem]{Definition}
\newtheorem{proposition}[theorem]{Proposition}
\newtheorem{remark}[theorem]{Remark}
\newtheorem{example}[theorem]{Example}
\newtheorem{corollary}[theorem]{Corollary}

\begin{document}

\title{Strictness of the containment of Kirk's asymptotic contractions in a class defined by two quasi-metrics}

\author{%
  Jie Shi\corrauth
}
\shortauthors{the Author(s)}

\address{%
 {Department of Mathematics and Statistics, Hubei Engineering University,
 272 Traffic Avenue, Xiaogan, 432000, Hubei, China}
}
\corraddr{Email: shijie@whu.edu.cn.
}

\begin{abstract}
Kirk's asymptotic contractions are defined by moduli $\psi_n$ which
converge uniformly on $[0,\infty)$ to a Boyd--Wong function and are
evaluated at the distance $d(x,y)$ only. A relaxed class, in which the
estimate involves a lower and an upper quasi-metric built from the four
distances $d(x,y)$, $d(x,Ty)$, $d(Tx,y)$ and $d(Tx,Ty)$, was recently
shown to contain Kirk's class and to admit the same fixed point
conclusion. It was left open whether the containment is strict. We prove
that it is. The separating maps form a family of radial damping maps
$T(x)=\lambda(\norm{x})x$ on a normed space, where
$\lambda\colon[0,\infty)\to[0,\rho]$ is continuous, $0\le\rho<1$,
nonconstant and periodic: every such map is a relaxed asymptotic
contraction with explicit linear moduli, whereas its one-step dispersion
at a fixed distance is unbounded, so that no finite Kirk modulus exists.
A slowly modulated nonperiodic damping law is treated as well. We further
show that strictness can occur only when the image $T(X)$ is not
relatively compact: in a complete space, if $T(X)$ is relatively compact,
the two classes coincide, hence on a compact (in particular finite) space
no separating example exists. Applications to semi-active vibration
isolation with amplitude-dependent damping and to digital maturity gaps in the digital transformation of small and medium-sized enterprises are discussed, and numerical illustrations are provided.
\end{abstract}


\maketitle

\section{Introduction}\label{sec:intro}

The contraction mapping principle of Banach \cite{Banach1922} requires a
uniform Lipschitz constant. A decisive generalization was obtained by
Boyd and Wong \cite{BoydWong1969}, who replaced the linear bound by a
nonlinear comparison function: if
\begin{equation}\label{eq:bw}
d(Tx,Ty)\le\psi\bigl(d(x,y)\bigr)
\end{equation}
for all $x,y$, where $\psi\colon[0,\infty)\to[0,\infty)$ is upper
semicontinuous, $\psi(0)=0$ and $\psi(t)<t$ for every $t>0$, then $T$ has
a unique fixed point and the Picard iterates converge to it.

Kirk \cite{Kirk2003} introduced the asymptotic version of \eqref{eq:bw}.
A continuous map $T$ of a metric space $(X,d)$ is an \emph{asymptotic
contraction} if there are functions
$\psi_n\colon[0,\infty)\to[0,\infty)$ such that
\begin{equation}\label{eq:kirk}
d(T^nx,T^ny)\le\psi_n\bigl(d(x,y)\bigr)
\qquad\text{for all }x,y\in X\text{ and }n\in\N,
\end{equation}
where $(\psi_n)$ converges to a Boyd--Wong function $\psi$ uniformly on
$[0,\infty)$. Kirk proved that if some orbit of $T$ is bounded, then $T$
has a unique fixed point $z$, and $T^nx\to z$ for every $x\in X$.

Numerous refinements followed; see
\cite{Suzuki2005,KirkXu2008,LindstromRoss2023,Arandelovic2005} and the
references therein. In all of them the estimate is a bound of
$d(T^nx,T^ny)$ in terms of the gap $d(x,y)$: a modulus of the gap in the
uniform case \eqref{eq:kirk}, and a state-dependent multiple of the gap in
the pointwise variants of \cite{KirkXu2008}.

In \cite{Shi2026} a relaxed class was proposed, in which the estimate
involves two quantities built directly from the map: a lower quasi-metric
$L_T$ and an upper quasi-metric $U_T$ which contains, besides the gap
$d(x,y)$, the one-step displacements $d(Tx,x)$ and $d(Ty,y)$. For a
self-map $T$ of a metric space $(X,d)$ put
\begin{equation}\label{eq:P}
P_T(x,y):=\min\bigl\{d(x,y),\,d(x,Ty),\,d(Tx,y),\,d(Tx,Ty)\bigr\},
\end{equation}
and define the \emph{lower} and \emph{upper quasi-metrics}
\begin{equation}\label{eq:L}
L_T(x,y):=
\begin{cases}
P_T(x,y), & \text{if }P_T(x,y)>0,\\[2pt]
d(Tx,Ty), & \text{if }P_T(x,y)=0,
\end{cases}
\end{equation}
\begin{equation}\label{eq:U}
U_T(x,y):=\max\bigl\{d(x,y),\,d(Tx,x),\,d(Ty,y),\,d(Tx,Ty)\bigr\}.
\end{equation}
One always has
\begin{equation}\label{eq:basic}
L_T(x,y)\le d(Tx,Ty)\le U_T(x,y),
\qquad d(x,y)\le U_T(x,y).
\end{equation}
The map $T$ is called a \emph{relaxed asymptotic contraction}
\cite{Shi2026} if there are functions
$\psi_n\colon[0,\infty)\to[0,\infty)$ such that
\begin{equation}\label{eq:relaxed}
L_T(T^nx,T^ny)\le\psi_n\bigl(U_T(x,y)\bigr)
\qquad\text{for all }x,y\in X\text{ and }n\in\N,
\end{equation}
and $(\psi_n)$ converges uniformly on every bounded subset of
$[0,\infty)$ to a Boyd--Wong function $\psi$. By \eqref{eq:basic},
condition \eqref{eq:relaxed} is weaker than \eqref{eq:kirk}: if
\eqref{eq:kirk} holds with moduli $\psi_n$, then, after replacing
$\psi_n$ by its monotone majorant, \eqref{eq:relaxed} holds with the same
moduli (see Remark~\ref{rem:monotone}). The main theorem of
\cite{Shi2026} states that the conclusion of Kirk's theorem persists for
this weaker class: a relaxed asymptotic contraction of a complete metric
space which has a bounded orbit admits a unique fixed point, to which all
iterates converge. The question whether the containment of Kirk's class
in the relaxed class is strict was left open in \cite{Shi2026}.

The present paper answers this question. Our results are the following.
\begin{enumerate}
\item Every continuous radial damping map $T(x)=\lambda(\norm{x})x$ on a
normed space, with $0\le\lambda\le\rho<1$, is a relaxed asymptotic
contraction, with the explicit linear moduli
$\psi_n(t)=2\rho^{\,n+1}(1-\rho)^{-1}t$ (Theorem~\ref{thm:family}).
\item If $\lambda$ is nonconstant and periodic, then $T$ is \emph{not} a
Kirk asymptotic contraction (Theorem~\ref{thm:dispersion} and
Corollary~\ref{cor:periodic}); hence the containment is strict
(Theorem~\ref{thm:strict}). Two variants, a slowly modulated nonperiodic
law and a two-dimensional law, are treated in Examples~\ref{ex:slow} and
\ref{ex:plane}.
\item Let $X$ be complete. If $T(X)$ is relatively compact, then every
relaxed asymptotic contraction is a Kirk asymptotic contraction
(Theorem~\ref{thm:compact}); hence on compact and on finite metric spaces
the two classes coincide, and strictness is necessarily a non-proper
phenomenon.
\item Two application scenarios are discussed in
Section~\ref{sec:applic}: semi-active vibration isolation with
amplitude-dependent damping, and digital maturity gaps in the digital transformation of small and medium-sized enterprises.
Numerical illustrations are given in Section~\ref{sec:applic} as well.
\end{enumerate}
The paper is organized as follows.
Section~\ref{sec:prelim} collects the definitions and the elementary
properties. Section~\ref{sec:family} introduces the radial damping maps,
proves that they belong to the relaxed class, and shows that they fail
Kirk's condition; this yields the strict containment.
Section~\ref{sec:compact} treats relatively compact images and compact
spaces, where the two classes coincide. Section~\ref{sec:applic} presents
the applications, and Section~\ref{sec:discussion} contains the conclusion
and some open problems.

\medskip
\noindent\textbf{Notation.} Throughout, $(X,d)$ is a metric space,
$\norm{\cdot}$ is the norm of a real normed space $X$ and
$d(x,y)=\norm{x-y}$; $e$ denotes a unit vector. For a map $T$ and
$x\in X$ we write $x_n:=T^nx$ and $\mathcal O(x):=\{T^nx:n\ge0\}$. All
metric notions are understood with respect to $d$.

\section{Preliminaries}\label{sec:prelim}

\subsection{The two quasi-metrics and the two classes}

\begin{definition}[Four-point minimum]\label{def:P}
For $x,y\in X$ define $P_T(x,y)$ by \eqref{eq:P}.
\end{definition}

\begin{definition}[Lower and upper quasi-metrics]\label{def:LU}
For $x,y\in X$ define $L_T(x,y)$ and $U_T(x,y)$ by \eqref{eq:L} and
\eqref{eq:U}, respectively.
\end{definition}

The terminology is justified by the following elementary properties.

\begin{lemma}\label{lem:basic}
For all $x,y\in X$:
\begin{enumerate}
\item $P_T(x,y)\ge0$ and $P_T(x,y)=P_T(y,x)$;
\item $L_T(x,y)\ge0$, $L_T(x,y)=L_T(y,x)$ and
$L_T(x,y)\le d(Tx,Ty)\le U_T(x,y)$;
\item $U_T(x,y)\ge0$, $U_T(x,y)=U_T(y,x)$ and $d(x,y)\le U_T(x,y)$.
\end{enumerate}
\end{lemma}

\begin{proof}
(i) The minimum of four nonnegative numbers is nonnegative. Each of the
four distances in \eqref{eq:P} is invariant under the interchange of $x$
and $y$: this is clear for $d(x,y)$ and $d(Tx,Ty)$, and the two remaining
terms are interchanged. Hence $P_T$ is symmetric.

(ii) If $P_T(x,y)>0$, then $L_T(x,y)=P_T(x,y)$ is one of the four
distances in \eqref{eq:P}, in particular $L_T(x,y)\le d(Tx,Ty)$. If
$P_T(x,y)=0$, then $L_T(x,y)=d(Tx,Ty)$ by definition. Thus
$L_T(x,y)\le d(Tx,Ty)$ in both cases; the symmetry follows from (i) and
the symmetry of $d(Tx,Ty)$; nonnegativity is clear.

(iii) $U_T(x,y)$ is a maximum of four nonnegative numbers, hence it is
nonnegative, symmetric, and dominates $d(x,y)$ and $d(Tx,Ty)$.
\end{proof}

The following estimate is the only place where the switching rule in the
definition of $L_T$ is used.

\begin{lemma}[Safe estimate]\label{lem:safe}
Let $u,v\in X$ be such that $P_T(u,v)>0$. Then
\begin{equation}\label{eq:safe}
d(Tu,Tv)\le P_T(u,v)+d(Tu,u)+d(Tv,v).
\end{equation}
\end{lemma}

\begin{proof}
Since $P_T(u,v)$ is the minimum of the four distances in \eqref{eq:P}, it
equals one of them; we distinguish the four cases.
If $P_T(u,v)=d(u,v)$, then the triangle inequality gives
$d(Tu,Tv)\le d(u,v)+d(Tu,u)+d(Tv,v)$, which is \eqref{eq:safe}.
If $P_T(u,v)=d(u,Tv)$, then
$d(Tu,Tv)\le d(u,Tv)+d(Tu,u)\le P_T(u,v)+d(Tu,u)+d(Tv,v)$.
If $P_T(u,v)=d(Tu,v)$, then
$d(Tu,Tv)\le d(Tu,v)+d(Tv,v)\le P_T(u,v)+d(Tu,u)+d(Tv,v)$.
Finally, if $P_T(u,v)=d(Tu,Tv)$, then
$d(Tu,Tv)=P_T(u,v)\le P_T(u,v)+d(Tu,u)+d(Tv,v)$.
This proves \eqref{eq:safe}.
\end{proof}

\begin{definition}[Relaxed asymptotic contraction]\label{def:relaxed}
A continuous map $T\colon X\to X$ is a \emph{relaxed asymptotic
contraction} if there are functions
$\psi_n\colon[0,\infty)\to[0,\infty)$ such that \eqref{eq:relaxed} holds
and $(\psi_n)$ converges to a Boyd--Wong function $\psi$ uniformly on
every bounded subset of $[0,\infty)$.
\end{definition}

\begin{definition}[Kirk asymptotic contraction]\label{def:kirk}
A continuous map $T\colon X\to X$ is a \emph{Kirk asymptotic contraction}
if there are functions $\psi_n\colon[0,\infty)\to[0,\infty)$ such that
\eqref{eq:kirk} holds and $(\psi_n)$ converges to a Boyd--Wong function
$\psi$ uniformly on $[0,\infty)$.
\end{definition}

\begin{remark}[Monotone moduli]\label{rem:monotone}
In Definitions~\ref{def:relaxed} and~\ref{def:kirk} the moduli may always
be assumed nondecreasing. Indeed, given $(\psi_n)$ satisfying the
corresponding inequality, put
\begin{equation}\label{eq:majorant}
h_n(t):=\sup_{0\le s\le t}\psi_n(s),\qquad t\ge0 .
\end{equation}
Then $h_n$ is nondecreasing and $h_n\ge\psi_n$, so that the corresponding
inequality holds with $h_n$ in place of $\psi_n$. Moreover, if $\psi_n\to\psi$
uniformly on $[0,M]$, then for $0\le t\le M$,
\begin{equation}\label{eq:majorant_conv}
\bigl|h_n(t)-g(t)\bigr|
\le\sup_{0\le s\le M}\bigl|\psi_n(s)-\psi(s)\bigr|,
\end{equation}
where $g(t)=\sup_{0\le s\le t}\psi(s)$; hence $h_n\to g$ uniformly on
$[0,M]$. By the argument of \cite{Shi2026}, $g$ is
nondecreasing and right continuous, $g(0)=0$ and $g(t)<t$ for $t>0$;
thus $g$ is a Boyd--Wong function. Consequently, replacing $\psi_n$ by
$h_n$ if necessary, we may assume the moduli in Definition~\ref{def:kirk}
to be nondecreasing.

With this normalization, \eqref{eq:kirk} implies \eqref{eq:relaxed}: by
Lemma~\ref{lem:basic},
\begin{equation}\label{eq:containment}
L_T(T^nx,T^ny)\le d(T^{n+1}x,T^{n+1}y)
\le\psi_{n+1}\bigl(d(x,y)\bigr)\le\psi_{n+1}\bigl(U_T(x,y)\bigr),
\end{equation}
and $\psi_{n+1}\to g$ uniformly on bounded intervals. Thus every Kirk
asymptotic contraction is a relaxed asymptotic contraction.
\end{remark}

\begin{remark}[Local versus global uniformity]\label{rem:local}
The requirement of uniformity on $[0,\infty)$ in \eqref{eq:kirk} is not
essential for the fixed point theorem. Indeed, by
Remark~\ref{rem:monotone} we may assume the moduli to be nondecreasing.
Let $x_0$ have a bounded orbit, put $x_n=T^nx_0$ and
$a_n=\sup_{p\ge1}d(x_n,x_{n+p})$, and let $D\ge a_n$ for all $n$. By
\eqref{eq:kirk} and the monotonicity of $\psi_m$,
\[
d(T^mx_{n+1},T^mx_{n+p+1})
\le\psi_m\bigl(d(x_{n+1},x_{n+p+1})\bigr)
\le\psi_m\bigl(a_{n+1}\bigr)
\le g\bigl(a_{n+1}\bigr)+\varepsilon
\]
for $m\ge M$, where $g$ is the monotone majorant of $\psi$ from
Remark~\ref{rem:monotone} and the last inequality follows from the
uniform convergence $\psi_m\to\psi$ on $[0,D]$ and from $\psi\le g$.
Hence $a_{n+m+1}\le g(a_{n+1})+\varepsilon$. Taking the limit
superior as $n\to\infty$ and using the right continuity of $g$ gives
$L\le g(L)+\varepsilon$, where $L=\limsup_na_n$; letting
$\varepsilon\to0^+$ and using $g(L)<L$ for $L>0$, we get $L=0$. Thus the
conclusion of Kirk's theorem holds under local uniformity as well. In
particular, the strictness established below is not
an artefact of the weaker mode of convergence used in
Definition~\ref{def:relaxed}.
\end{remark}

\subsection{The fixed point theorem for the relaxed class}

The following result is the main theorem of \cite{Shi2026}; we quote it
and use it only in Section~\ref{sec:compact}.

\begin{theorem}[\cite{Shi2026}]\label{thm:quoted}
Let $(X,d)$ be a complete metric space and let $T\colon X\to X$ be a
relaxed asymptotic contraction. If some orbit of $T$ is bounded, then $T$
has a unique fixed point $z$, and $T^nx\to z$ for every $x\in X$.
\end{theorem}

\begin{remark}\label{rem:quoted}
For the radial maps of Section~\ref{sec:family} the conclusion of
Theorem~\ref{thm:quoted} can also be verified directly, see
Corollary~\ref{cor:radial_convergence}; the point of
Theorem~\ref{thm:quoted} is that it applies to arbitrary relaxed
asymptotic contractions and is needed in Section~\ref{sec:compact}.
\end{remark}

\section{Radial damping maps}\label{sec:family}

Radial damping maps play a central role in this paper. On the one hand,
they provide an explicit model of the relaxed class: whenever the damping
law is continuous and bounded by a constant $\rho<1$, the associated map
is a relaxed asymptotic contraction with explicit linear moduli
(Theorem~\ref{thm:family}). On the other hand, when the damping law is
nonconstant and periodic, the map is not a Kirk asymptotic contraction
(Theorem~\ref{thm:dispersion} and Corollary~\ref{cor:periodic}): its
one-step dispersion at a fixed distance is unbounded. The radial family
therefore supplies the explicit witnesses of the strict containment
(Theorem~\ref{thm:strict}), and the analysis of this section is the core
of the paper.

\subsection{Definition and main structural result}

\begin{definition}[Damping law and radial damping map]\label{def:family}
A \emph{damping law} is a continuous function
$\lambda\colon[0,\infty)\to[0,\rho]$ with $0\le\rho<1$. The associated
\emph{radial damping map} $T=T_\lambda\colon X\to X$ on a real normed
space $X$ is defined by
\begin{equation}\label{eq:T}
T(x):=\lambda\bigl(\norm{x}\bigr)\,x,\qquad x\in X .
\end{equation}
\end{definition}

Physically, \eqref{eq:T} is the one-step map of a one-degree-of-freedom
oscillator in which the damping ratio, over one sampling period, depends
on the amplitude of the displacement. The assumption
$0\le\lambda\le\rho<1$ means that the damping is never negative and that
the map is never expansive in the radial direction; the state dependence
of $\lambda$ is the feature which distinguishes the present situation
from a linear contraction, see Section~\ref{sec:applic}.

\begin{theorem}[Radial damping maps]\label{thm:family}
Let $X$ be a real normed space, let $\lambda\colon[0,\infty)\to[0,\rho]$
be continuous with $0\le\rho<1$, and let $T=T_\lambda$ be given by
\eqref{eq:T}. Then:
\begin{enumerate}
\item $T$ is continuous and
\begin{equation}\label{eq:radial_bounds}
\norm{T(x)}\le\rho\norm{x},
\qquad
\norm{T(x)-x}\ge(1-\rho)\norm{x}
\qquad(x\in X);
\end{equation}
\item every orbit of $T$ converges to $0$, and $0$ is the unique fixed
point of $T$;
\item $T$ is a relaxed asymptotic contraction in the sense of
Definition~\ref{def:relaxed}. More precisely,
\begin{equation}\label{eq:family_est}
L_T(T^nx,T^ny)\le\psi_n\bigl(U_T(x,y)\bigr)
\qquad(x\in X,\ y\in X,\ n\ge0),
\end{equation}
where
\begin{equation}\label{eq:psi}
\psi_n(t):=\frac{2\rho^{\,n+1}}{1-\rho}\,t,
\qquad t\ge0 ;
\end{equation}
the sequence $(\psi_n)$ converges to $\psi\equiv0$ uniformly on every
bounded subset of $[0,\infty)$, and $\psi$ satisfies the Boyd--Wong
condition.
\end{enumerate}
\end{theorem}

\begin{proof}
\noindent\textbf{Step 1: continuity and the two bounds.}
The function $x\mapsto\norm{x}$ is continuous, hence so is
$x\mapsto\lambda(\norm{x})$; being a product of continuous maps, $T$ is
continuous. Since $0\le\lambda(\norm{x})\le\rho$, we have
$\norm{T(x)}=\lambda(\norm{x})\norm{x}\le\rho\norm{x}$. Moreover,
\[
\norm{T(x)-x}
=\bigl|\lambda(\norm{x})-1\bigr|\,\norm{x}
=\bigl(1-\lambda(\norm{x})\bigr)\norm{x}
\ge(1-\rho)\norm{x},
\]
because $\lambda(\norm{x})\le\rho<1$.

\noindent\textbf{Step 2: orbits and fixed points.}
By Step 1, $\norm{T^nx}\le\rho^n\norm{x}$ for all $n\ge0$, hence
$T^nx\to0$. If $T(x)=x$, then
$\bigl(\lambda(\norm{x})-1\bigr)x=0$; since
$\lambda(\norm{x})-1\le\rho-1<0$, this forces $x=0$. Conversely $T(0)=0$.
Hence $0$ is the unique fixed point.

\noindent\textbf{Step 3: a lower bound for $U_T$.}
By \eqref{eq:radial_bounds} and the definition of $U_T$,
\[
U_T(x,y)\ge d(Tx,x)=\norm{T(x)-x}\ge(1-\rho)\norm{x},
\qquad
U_T(x,y)\ge d(Ty,y)\ge(1-\rho)\norm{y}.
\]
Therefore
\begin{equation}\label{eq:sum_bound}
\norm{x}+\norm{y}
\le\frac{2}{1-\rho}\,U_T(x,y)
\qquad(x,y\in X).
\end{equation}

\noindent\textbf{Step 4: verification of \eqref{eq:family_est}.}
Fix $x,y\in X$ and $n\ge0$. By Lemma~\ref{lem:basic},
$L_T(u,v)\le d(Tu,Tv)$ for all $u,v\in X$; applying this with
$u=T^nx$ and $v=T^ny$, and then using Step 1 and
\eqref{eq:sum_bound}, we obtain
\begin{equation}\label{eq:family_chain}
\begin{aligned}
L_T(T^nx,T^ny)
&\le d(T^{n+1}x,T^{n+1}y)
=\norm{T^{n+1}x-T^{n+1}y}\\
&\le\norm{T^{n+1}x}+\norm{T^{n+1}y}
\le\rho^{\,n+1}\bigl(\norm{x}+\norm{y}\bigr)
\le\frac{2\rho^{\,n+1}}{1-\rho}\,U_T(x,y),
\end{aligned}
\end{equation}
which is \eqref{eq:family_est} with $\psi_n$ given by \eqref{eq:psi}.

\noindent\textbf{Step 5: the moduli.}
Each $\psi_n$ is nonnegative and $\psi_n(0)=0$. If $t\le M$, then
$0\le\psi_n(t)=2\rho^{n+1}(1-\rho)^{-1}t\le2M(1-\rho)^{-1}\rho^{n+1}\to0$,
so $(\psi_n)$ converges to $\psi\equiv0$ uniformly on $[0,M]$. The
function $\psi\equiv0$ is continuous, $\psi(0)=0$ and $\psi(t)=0<t$ for
$t>0$, so $\psi$ satisfies the Boyd--Wong condition. This completes the
proof.
\end{proof}

\begin{corollary}[Global convergence of radial damping maps]
\label{cor:radial_convergence}
Under the hypotheses of Theorem~\ref{thm:family}, for every $x_0\in X$,
\begin{equation}\label{eq:rate}
\norm{T^nx_0}\le\rho^n\norm{x_0}
\qquad(n\ge0),
\end{equation}
so that $T^nx_0\to0$. In particular, the equilibrium $0$ of the damping
law \eqref{eq:T} is globally asymptotically stable, with geometric rate
$\rho$.
\end{corollary}

\begin{proof}
The estimate \eqref{eq:rate} follows by induction from the first bound in
\eqref{eq:radial_bounds} (Step 1 of the proof of
Theorem~\ref{thm:family}); the uniqueness of the equilibrium is Step 2
there.
\end{proof}

\begin{remark}[The role of the explicit moduli]\label{rem:role}
In \eqref{eq:family_chain} the argument of the modulus is
$U_T(x,y)$, which dominates $(1-\rho)(\norm{x}+\norm{y})/2$. Thus the
relaxed condition compares $L_T(T^nx,T^ny)$ with a quantity of the same
order of magnitude as the states: the upper quasi-metric carries the
information about \emph{where} the pair sits in the space, not only about
the gap between its points. In contrast, Kirk's condition
\eqref{eq:kirk} evaluates the modulus at $d(x,y)$ alone. This is the
structural reason for the strictness proved in
this section.
\end{remark}

\begin{example}[Constant damping law]\label{ex:constant}
Let $\lambda\equiv c$ with $0\le c<1$. Then $T(x)=cx$, and \eqref{eq:T}
is a linear contraction with constant $c$. Consequently $T$ is a Kirk
asymptotic contraction: with $\psi_0(t)=t$ and $\psi_n(t)=ct$ for
$n\ge1$, the inequality \eqref{eq:kirk} holds because $c^n\le c$, and
$(\psi_n)$ converges uniformly on $[0,\infty)$ to the Boyd--Wong function
$\psi(t)=ct$. Hence $T$ is also a relaxed one. This shows that the radial
family of Theorem~\ref{thm:family} contains the classical linear
contractions, for which the two classes coincide.
\end{example}

\begin{figure}[htbp]
\centering
\includegraphics[width=0.88\textwidth]{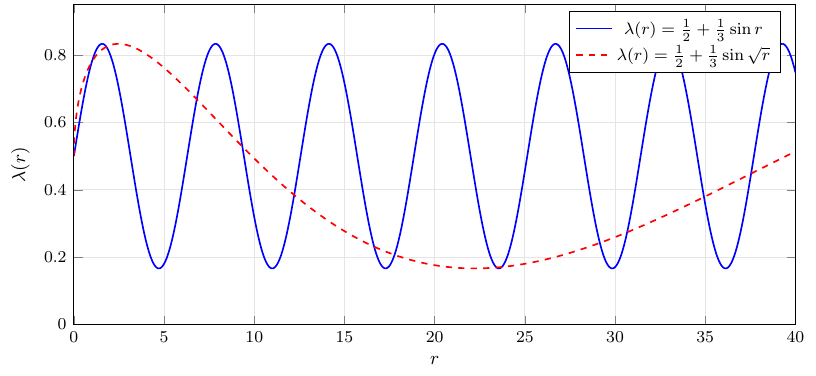}
\caption{Two damping laws of the family \eqref{eq:T}. The first one is
periodic and nonconstant (Example~\ref{ex:harmonic}); the second one is
nonperiodic and slowly modulated. Both are relaxed asymptotic
contractions.}
\label{fig:laws}
\end{figure}

\begin{example}[Harmonic damping law]\label{ex:harmonic}
Let
\begin{equation}\label{eq:harmonic}
\lambda(r):=\frac12+\frac13\sin r,\qquad r\ge0 .
\end{equation}
Then $\lambda$ is continuous, $2\pi$-periodic, nonconstant, and
\begin{equation}\label{eq:harmonic_range}
0\le\lambda(r)\le\frac56=:\rho,\qquad r\ge0 .
\end{equation}
By Theorem~\ref{thm:family}, $T(x)=\bigl(\frac12+\frac13\sin\norm{x}\bigr)x$
is a relaxed asymptotic contraction with
\begin{equation}\label{eq:harmonic_moduli}
\psi_n(t)=12\Bigl(\frac56\Bigr)^{n+1}t,
\end{equation}
and every orbit converges to $0$ geometrically. The map is not a linear
contraction: its one-dimensional derivative is
\begin{equation}\label{eq:harmonic_derivative}
T'(x)=\frac12+\frac13\sin x+\frac{x}{3}\cos x,
\end{equation}
which is unbounded on $[0,\infty)$. Below we show that $T$ is not a Kirk asymptotic
contraction; the profiles of this law, together with the slowly modulated
law treated below, are plotted in Figure~\ref{fig:laws}.
\end{example}

\subsection{Failure of Kirk's condition}\label{sec:strict}

The obstruction to Kirk's condition is the behaviour of the map on pairs
at a \emph{fixed} distance: Kirk's inequality with $n=1$ requires the
quantity
\begin{equation}\label{eq:dispersion}
\Delta_T(t):=\sup\bigl\{d(Tx,Ty): d(x,y)=t\bigr\},
\qquad t>0,
\end{equation}
to be finite, and more precisely it forces
$\Delta_T(t)\le\varphi_1(t)$ for every $t$. This motivates the following
criterion.

\begin{theorem}[Unbounded dispersion]\label{thm:dispersion}
Let $X$ be a real normed space with $\dim X\ge1$, let
$\lambda\colon[0,\infty)\to[0,\rho]$ be continuous with $0\le\rho<1$, and
let $T$ be given by \eqref{eq:T}. Suppose that for some $t>0$,
\begin{equation}\label{eq:dispersion_inf}
\sup_{r\ge0}
\bigl|r\lambda(r)-(r+t)\lambda(r+t)\bigr|=+\infty .
\end{equation}
Then $T$ is not a Kirk asymptotic contraction: there are no functions
$\varphi_n\colon[0,\infty)\to[0,\infty)$ such that \eqref{eq:kirk} holds.
\end{theorem}

\begin{proof}
Assume, for contradiction, that such functions $\varphi_n$ exist, and
let $e\in X$ be a unit vector. For every $r\ge0$ put
\begin{equation}\label{eq:pair}
x:=re,\qquad y:=(r+t)e .
\end{equation}
Then $d(x,y)=\norm{t e}=t$, while, by \eqref{eq:T},
\[
T(x)=\lambda(r)\,re=\lambda(r)\,x,
\qquad
T(y)=\lambda(r+t)\,(r+t)e=\lambda(r+t)\,y,
\]
and therefore
\begin{equation}\label{eq:dispersion_pair}
d(Tx,Ty)
=\bigl\|\lambda(r)\,re-\lambda(r+t)\,(r+t)e\bigr\|
=\bigl|r\lambda(r)-(r+t)\lambda(r+t)\bigr| .
\end{equation}
Applying \eqref{eq:kirk} with $n=1$ to the pair \eqref{eq:pair} gives
\begin{equation}\label{eq:dispersion_bound}
\bigl|r\lambda(r)-(r+t)\lambda(r+t)\bigr|
\le\varphi_1(t)
\qquad\text{for every }r\ge0 .
\end{equation}
Since $\varphi_1(t)$ is a finite real number, this contradicts
\eqref{eq:dispersion_inf}. Hence no family $(\varphi_n)$ satisfying
\eqref{eq:kirk} exists.
\end{proof}

\begin{corollary}[Periodic damping laws]\label{cor:periodic}
Let $X$ be a real normed space with $\dim X\ge1$, let
$\lambda\colon[0,\infty)\to[0,\rho]$ be continuous, $p$-periodic and
nonconstant, with $0\le\rho<1$, and let $T$ be the associated radial
damping map \eqref{eq:T} on $X$. Then $T$ is a relaxed asymptotic contraction
but not a Kirk asymptotic contraction.
\end{corollary}

\begin{proof}
By Theorem~\ref{thm:family}, $T$ is a relaxed asymptotic contraction.

We show that \eqref{eq:dispersion_inf} holds for a suitable $t>0$. Since
$\lambda$ is continuous and nonconstant, there exists
$t\in(0,p)$ which is not a period of $\lambda$. Indeed, if every
$t\in(0,p)$ were a period, then for arbitrary $0\le r<s$ choose
$N\in\N$ with $(s-r)/N<p$; then $\lambda(r)=\lambda(s)$, because
$s-r=N\cdot\frac{s-r}{N}$ is a multiple of the period
$\frac{s-r}{N}$, and $\lambda$ would be constant, a contradiction.

Fix such a $t$ and put
\begin{equation}\label{eq:g}
g(r):=\lambda(r+t)-\lambda(r),\qquad r\ge0 .
\end{equation}
Then $g$ is continuous, $p$-periodic, and not identically zero. Choose
$r_0\ge0$ with $c:=g(r_0)\ne0$, and for $k\ge1$ put $r_k:=r_0+kp$. By
periodicity, $g(r_k)=c$ for every $k$. Hence
\begin{equation}\label{eq:dispersion_along}
\begin{aligned}
\bigl|r_k\lambda(r_k)-(r_k+t)\lambda(r_k+t)\bigr|
&=\bigl|-r_k\bigl(\lambda(r_k+t)-\lambda(r_k)\bigr)
-t\lambda(r_k+t)\bigr| \\
&=\bigl|-r_k\,c-t\,\lambda(r_k+t)\bigr|
\ge r_k|c|-t\rho,
\end{aligned}
\end{equation}
which tends to $+\infty$ as $k\to\infty$. Thus \eqref{eq:dispersion_inf}
holds, and Theorem~\ref{thm:dispersion} shows that $T$ is not a Kirk
asymptotic contraction.
\end{proof}

\begin{remark}\label{rem:criterion}
The criterion of Theorem~\ref{thm:dispersion} concerns only the values of
$T$ on the space, not the iterates: the obstruction appears already at
$n=1$. For this reason the conclusion of Corollary~\ref{cor:periodic}
remains true if in \eqref{eq:kirk} the uniform convergence on
$[0,\infty)$ is replaced by local uniform convergence, or dropped
altogether. In particular, the strictness proved below is not related to
the local uniformity discussed in Remark~\ref{rem:local}.
\end{remark}

The following examples make the failure explicit.

\begin{figure}[htbp]
\centering
\includegraphics[width=0.88\textwidth]{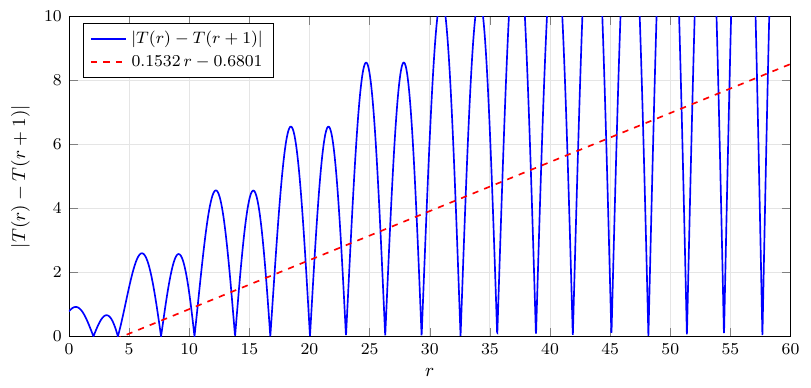}
\caption{One-step dispersion at distance $1$ for the harmonic law
\eqref{eq:harmonic}. The dispersion is unbounded, so that Kirk's modulus
would have to take the value $+\infty$ at $t=1$.}
\label{fig:disp}
\end{figure}

\begin{example}[Harmonic damping: explicit divergence]\label{ex:sin}
Let $X=\R$ and let $T$ be given by the harmonic law
\eqref{eq:harmonic}. By Example~\ref{ex:harmonic}, $T$ is a relaxed
asymptotic contraction. Here we exhibit explicitly the divergence
required by Theorem~\ref{thm:dispersion}, and we compute some numerical
values.

For $k\ge0$ put
\begin{equation}\label{eq:xk}
p_k:=\frac{\pi}{2}+2k\pi,\qquad q_k:=p_k+1 .
\end{equation}
Then $d(p_k,q_k)=1$, $\sin p_k=1$ and $\sin q_k=\sin(p_k+1)=\cos1$.
Hence
\begin{equation}\label{eq:sin_disp}
\begin{aligned}
T(p_k)-T(q_k)
&=\Bigl(\frac12+\frac13\Bigr)p_k
-\Bigl(\frac12+\frac13\cos1\Bigr)(p_k+1) \\
&=\frac{1-\cos1}{3}\,p_k-\Bigl(\frac12+\frac{\cos1}{3}\Bigr),
\end{aligned}
\end{equation}
and since $1-\cos1>0$, this tends to $+\infty$ as $k\to\infty$. Thus
$\Delta_T(1)=+\infty$, and Theorem~\ref{thm:dispersion} (with $t=1$)
shows that $T$ is not a Kirk asymptotic contraction.

Table~\ref{tab:disp} collects the numerical values of
$|T(p_k)-T(q_k)|$ for several $k$, and Figure~\ref{fig:disp} displays
the function $r\mapsto|T(r)-T(r+1)|$. Note the asymptotic formula
\begin{equation}\label{eq:sin_asymptotic}
T(p_k)-T(q_k)
=\frac{1-\cos1}{3}\,p_k-\Bigl(\frac12+\frac{\cos1}{3}\Bigr)
\approx0.1532\,p_k-0.6801,
\end{equation}
so that the dispersion at the fixed distance $1$ grows linearly in the
position.
\end{example}

\begin{table}[htbp]
\centering
\begin{tabular}{cccc}
\hline
$k$ & $p_k$ & $q_k$ & $|T(p_k)-T(q_k)|$ \\
\hline
$1$     & $7.853982$  & $8.853982$  & $0.523385$ \\
$10$    & $64.402649$ & $65.402649$ & $9.188482$ \\
$100$   & $629.889327$& $630.889327$& $95.839456$ \\
$1000$  & $6284.756104$& $6285.756104$& $962.349196$ \\
\hline
\end{tabular}
\caption{The one-step dispersion at the fixed distance $1$ for the
harmonic law \eqref{eq:harmonic}: the values of
$|T(p_k)-T(q_k)|$ with $p_k=\pi/2+2k\pi$ and $q_k=p_k+1$. The dispersion
grows linearly in $p_k$, in accordance with \eqref{eq:sin_asymptotic}.}
\label{tab:disp}
\end{table}

\begin{example}[Slowly modulated nonperiodic damping]\label{ex:slow}
Let
\begin{equation}\label{eq:slow}
\lambda(r):=\frac12+\frac13\sin\sqrt{r},\qquad r\ge0 .
\end{equation}
Then $\lambda$ is continuous, nonperiodic, and
$0\le\lambda(r)\le5/6$ for all $r\ge0$, so, by
Theorem~\ref{thm:family}, $T(x)=\bigl(\frac12+\frac13\sin\sqrt{\norm{x}}\bigr)x$
is a relaxed asymptotic contraction with the moduli
$\psi_n(t)=12(5/6)^{n+1}t$. We show that $T$ is not a Kirk asymptotic
contraction, so that periodicity is not essential for the phenomenon.

Take $t=1$ and, for $k\ge1$,
\begin{equation}\label{eq:rk}
r_k:=(2k\pi)^2 .
\end{equation}
Then $\sin\sqrt{r_k}=\sin(2k\pi)=0$, hence $\lambda(r_k)=1/2$. Next,
since
\begin{equation}\label{eq:sqrt_expansion}
\sqrt{r_k+1}
=2k\pi\sqrt{1+\frac{1}{4k^2\pi^2}}
=2k\pi+\frac{1}{4k\pi}+O\bigl(k^{-3}\bigr),
\end{equation}
we have
\begin{equation}\label{eq:sin_expansion}
\sin\sqrt{r_k+1}
=\sin\Bigl(\frac{1}{4k\pi}+O(k^{-3})\Bigr)
=\frac{1}{4k\pi}+O\bigl(k^{-3}\bigr),
\end{equation}
and therefore
\begin{equation}\label{eq:slow_lambda}
\lambda(r_k+1)=\frac12+\frac{1}{12k\pi}+O\bigl(k^{-3}\bigr).
\end{equation}
Consequently
\begin{equation}\label{eq:slow_disp}
\begin{aligned}
r_k\lambda(r_k)-(r_k+1)\lambda(r_k+1)
&=\frac{r_k}{2}-\frac{r_k+1}{2}
-\frac{r_k+1}{12k\pi}+O\bigl(k^{-1}\bigr) \\
&=-\frac12-\Bigl(\frac{k\pi}{3}+\frac{1}{12k\pi}\Bigr)+O\bigl(k^{-1}\bigr)
\longrightarrow-\infty,
\end{aligned}
\end{equation}
where we used $(r_k+1)/(12k\pi)=k\pi/3+1/(12k\pi)$. Thus
\eqref{eq:dispersion_inf} holds with $t=1$, and
Theorem~\ref{thm:dispersion} shows that $T$ is not a Kirk asymptotic
contraction.
\end{example}

\begin{example}[A two-dimensional separating example]\label{ex:plane}
Let $X=\R^2$ with the Euclidean norm and let
\begin{equation}\label{eq:plane}
T(x):=\Bigl(\frac12+\frac13\sin\norm{x}\Bigr)x .
\end{equation}
By Theorem~\ref{thm:family}, $T$ is a relaxed asymptotic contraction with
$\psi_n(t)=12(5/6)^{n+1}t$, and all orbits converge to $0$. On the ray
$\{re_1:r\ge0\}$, with $e_1=(1,0)$, the computation of
Example~\ref{ex:sin} applies verbatim: for $p_k=(\pi/2+2k\pi)e_1$ and
$q_k=p_k+e_1$ we have $d(p_k,q_k)=1$, $\sin\norm{p_k}=1$ and
$\sin\norm{q_k}=\cos1$, hence
\begin{equation}\label{eq:plane_disp}
T(p_k)-T(q_k)
=\frac{1-\cos1}{3}\,p_k-\Bigl(\frac12+\frac{\cos1}{3}\Bigr)e_1 ,
\qquad
\bigl\|T(p_k)-T(q_k)\bigr\|\longrightarrow\infty ,
\end{equation}
and $T$ is not a Kirk asymptotic contraction. Thus the phenomenon is not
a one-dimensional accident; it occurs in every normed space of
dimension at least one, and the radial damping law can be imposed in any
direction.
\end{example}

The examples above show that the relaxed class contains maps outside
Kirk's class. The following theorem records this fact.

\begin{theorem}[Strict containment]\label{thm:strict}
The class of relaxed asymptotic contractions strictly contains the class
of Kirk asymptotic contractions.
\end{theorem}

\begin{proof}
By Remark~\ref{rem:monotone}, every Kirk asymptotic contraction is a
relaxed asymptotic contraction. By Example~\ref{ex:sin} (or
Example~\ref{ex:slow}, or Example~\ref{ex:plane}) there exists a relaxed
asymptotic contraction which is not a Kirk asymptotic contraction.
Therefore the containment is strict. This answers the question left open
in \cite{Shi2026}.
\end{proof}

\section{Relatively compact images and compact spaces}\label{sec:compact}

The separating examples of Subsection~\ref{sec:strict} all exploit the
behaviour of the map ``at infinity'': the dispersion is generated by pairs
of states of large norm. The following result shows that this is not an
accident: whenever the image of $T$ is relatively compact, the relaxed
condition collapses to Kirk's condition.

\begin{theorem}[Coincidence for relatively compact images]
\label{thm:compact}
Let $(X,d)$ be a complete metric space and let $T\colon X\to X$ be a
relaxed asymptotic contraction. Assume that $T(X)$ is relatively compact.
Then $T$ is a Kirk asymptotic contraction. More precisely, let $z$ be the
unique fixed point of $T$ and, with $K:=\overline{T(X)}$, let
\begin{equation}\label{eq:epsilon}
\varepsilon_m:=\sup_{u\in K}d\bigl(T^mu,z\bigr),\qquad m\ge0 .
\end{equation}
Then $\varepsilon_m\to0$ and
\begin{equation}\label{eq:kirk_compact}
d(T^nx,T^ny)\le\varphi_n\bigl(d(x,y)\bigr)
\qquad(x,y\in X,\ n\ge0),
\end{equation}
where $\varphi_0(t)=t$ and $\varphi_n(t)\equiv2\varepsilon_{n-1}$ for
$n\ge1$. In particular, $\varphi_n\to0$ uniformly on $[0,\infty)$, so that
\eqref{eq:kirk_compact} is Kirk's condition.
\end{theorem}

\medskip
\noindent\textit{Strategy of the proof.} The proof has three steps. In
Step~1 we show that the increments $d(T^nx,T^{n+1}x)$ tend to $0$
uniformly in $x$. In Step~2 we show that the tails $d(T^nx,T^{n+p}x)$
tend to $0$ uniformly in $x$ and $p$; here the safe estimate
(Lemma~\ref{lem:safe}) and the monotone moduli
(Remark~\ref{rem:monotone}) are used. In Step~3 we convert these uniform
estimates into Kirk's condition. Relative compactness provides the compact
set $K$ on which all the asymptotic constants are uniform.

\begin{proof}
By Remark~\ref{rem:monotone} we may assume the moduli $\psi_n$ to be
nondecreasing. Since $T(X)$ is relatively compact, it is bounded, and for
every $x\in X$ the orbit satisfies $\mathcal O(x)\subset\{x\}\cup T(X)$,
hence it is bounded; by Theorem~\ref{thm:quoted}, $T$ has a unique fixed
point $z$ and $T^mu\to z$ for every $u\in X$. Put $K:=\overline{T(X)}$ and
$D:=\diam K<\infty$.

\noindent\textbf{Step 1: uniform decay of the increments.}
For $x\in X$ and $n\ge0$ put $d_n(x):=d(T^nx,T^{n+1}x)$ and
$\delta_n:=\sup_{x\in X}d_n(x)$. We claim that $\delta_n\to0$.
For $u:=T^nx$ and $v:=T^{n+1}x$ we have
$d(Tu,v)=d(T^{n+1}x,T^{n+1}x)=0$, hence $P_T(u,v)=0$ and therefore
$L_T(u,v)=d(Tu,Tv)=d_{n+1}(x)$; moreover
$U_T(u,v)=\max\{d_n(x),d_{n+1}(x)\}$. Applying \eqref{eq:relaxed} to the
pair $(u,v)$ with shift $m$ gives
\begin{equation}\label{eq:inc}
d_{n+m+1}(x)\le\psi_m\bigl(\max\{d_n(x),d_{n+1}(x)\}\bigr)
\qquad(x\in X,\ n,m\ge0).
\end{equation}
Taking $m=n$ and $m=n-1$ ($n\ge1$), passing to the supremum over $x\in X$
and using the monotonicity of the moduli, we obtain
\begin{equation}\label{eq:inc2}
\delta_{2n+1}\le\psi_n\bigl(\max\{\delta_n,\delta_{n+1}\}\bigr),
\qquad
\delta_{2n}\le\psi_{n-1}\bigl(\max\{\delta_n,\delta_{n+1}\}\bigr).
\end{equation}
For $n\ge1$ one has $\delta_n\le D$, because $T^nx,T^{n+1}x\in T(X)\subset
K$. Let $R:=\limsup_{n\to\infty}\delta_n$. Passing to the limit superior in
\eqref{eq:inc2}, using the local uniform convergence of $\psi_m$ to $\psi$
on $[0,D]$ and the monotone majorant $g$ of $\psi$
(Remark~\ref{rem:monotone}), we get $R\le g(R)$; hence $R=0$, that is,
$\delta_n\to0$.

\noindent\textbf{Step 2: the uniform Cauchy property.}
For $x\in X$ and $n\ge0$ put
\begin{equation}\label{eq:an}
a_n(x):=\sup_{p\ge1}d(T^nx,T^{n+p}x),
\qquad
A_n:=\sup_{x\in X}a_n(x).
\end{equation}
We show that $A_n\to0$. Let $\varepsilon>0$. By Step 1 there is $N$ with
$d_k(x)\le\varepsilon$ for all $x\in X$ and all $k\ge N$. Fix $n\ge N$,
$p\ge1$ and $m\ge1$, and put $u:=T^{n+m}x$, $v:=T^{n+p+m}x$. If
$P_T(u,v)>0$, then Lemma~\ref{lem:safe} and \eqref{eq:relaxed} give
\[
d(Tu,Tv)\le P_T(u,v)+d(Tu,u)+d(Tv,v)
\le\psi_m\bigl(U_T(T^nx,T^{n+p}x)\bigr)+2\varepsilon,
\]
because $L_T(T^nx,T^{n+p}x)=P_T(T^nx,T^{n+p}x)\le
\psi_m(U_T(T^nx,T^{n+p}x))$ and $d(Tu,u)=d_{n+m}(x)\le\varepsilon$,
$d(Tv,v)=d_{n+p+m}(x)\le\varepsilon$. If $P_T(u,v)=0$, then
$L_T(u,v)=d(Tu,Tv)$ and the same inequality holds without the term
$P_T(u,v)$. Furthermore,
\[
U_T(T^nx,T^{n+p}x)\le a_n(x)+2\varepsilon,
\]
since its four entries are bounded by $a_n(x)$, $\varepsilon$,
$\varepsilon$ and $a_n(x)+2\varepsilon$. Consequently, for all $m\ge M$
with $M$ so large that $\psi_m\le\psi+\varepsilon\le g+\varepsilon$ on
$[0,3D]$,
\[
d(T^{n+m+1}x,T^{n+p+m+1}x)\le g\bigl(a_n(x)+2\varepsilon\bigr)+3\varepsilon.
\]
Taking the supremum over $p\ge1$ and then over $x\in X$, we obtain
\begin{equation}\label{eq:A_rec}
A_{n+m+1}\le g(A_n+2\varepsilon)+3\varepsilon
\qquad(n\ge N,\ m\ge M).
\end{equation}
Let $A:=\limsup_{n\to\infty}A_n$. Passing to the limit superior in
\eqref{eq:A_rec} and letting $\varepsilon\to0^+$ yields $A\le g(A)$; hence
$A=0$.

\noindent\textbf{Step 3: Kirk's condition.}
For $n\ge1$, $x\in X$ and $p\ge1$,
\[
d(T^nx,z)\le d(T^nx,T^{n+p}x)+d(T^{n+p}x,z)
\le A_n+d(T^{n+p}x,z);
\]
letting $p\to\infty$ gives $d(T^nx,z)\le A_n$. In particular
$\varepsilon_m\le A_m\to0$ for $m\ge1$. Moreover, for all $x,y\in X$ and
$n\ge1$,
\begin{equation}\label{eq:compact_est}
d(T^nx,T^ny)\le d(T^nx,z)+d(T^ny,z)\le2\varepsilon_{n-1}
=\varphi_n\bigl(d(x,y)\bigr),
\end{equation}
because $d(T^nx,z)=d(T^{n-1}(Tx),z)\le\varepsilon_{n-1}$ and $Tx\in K$.
For $n=0$ the inequality reads $d(x,y)\le\varphi_0(d(x,y))=d(x,y)$. Finally
$\varphi_n\to0$ uniformly on $[0,\infty)$, because
$\sup_{t\ge0}|\varphi_n(t)|=2\varepsilon_{n-1}\to0$, and $\varphi\equiv0$
is a Boyd--Wong function. Hence $T$ is a Kirk asymptotic contraction.
\end{proof}

\begin{corollary}[Compact spaces]\label{cor:compact}
If $X$ is compact, then the class of relaxed asymptotic contractions and
the class of Kirk asymptotic contractions coincide. In particular, no
separating example exists on a compact metric space, and hence none
exists on a finite metric space.
\end{corollary}

\begin{proof}
On a compact space the image $T(X)$ is relatively compact, so
Theorem~\ref{thm:compact} gives the inclusion of the relaxed class into
Kirk's class. The reverse inclusion is Remark~\ref{rem:monotone}.
\end{proof}

\begin{remark}[Images of radial damping maps]\label{rem:image}
Let $\lambda$ be a damping law with $\lambda(r)\ge c$ for some $c>0$ and
all $r\ge0$. Then
$\norm{T(x)}=\lambda(\norm{x})\norm{x}\ge c\norm{x}$, so the image $T(X)$
contains points of arbitrarily large norm and is therefore not relatively
compact. This is consistent with Theorem~\ref{thm:compact}: the
separating examples of Section~\ref{sec:strict} have unbounded (and hence
not relatively compact) images.

In a finite-dimensional space, a set is relatively compact if and only if
it is bounded. Hence, in $\R^n$, a separating example must have an
unbounded image; the harmonic law of Example~\ref{ex:sin} satisfies
$\norm{T(x)}=\lambda(\norm{x})\norm{x}\ge\frac16\norm{x}$, so its image
is all of $\R$ in the one-dimensional case. Combining this with
Theorem~\ref{thm:strict}, we conclude that the strictness of the
containment is a genuinely non-proper phenomenon: it is caused by the
unbounded dispersion of the map at infinity, not by any local pathology.
\end{remark}

\section{Applications}\label{sec:applic}

\subsection{Semi-active vibration isolation with amplitude-dependent
damping}\label{subsec:damper}

\begin{figure}[htbp]
\centering
\includegraphics[width=0.62\textwidth]{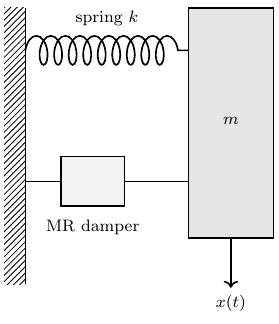}
\caption{A one-degree-of-freedom oscillator with a spring and a
magnetorheological damper. The damping ratio is modulated harmonically
with the displacement, which leads to the damping law
\eqref{eq:damper_law} and to the multiplicative sampled map
\eqref{eq:damper_map}.}
\label{fig:schematic}
\end{figure}

Consider a one-degree-of-freedom oscillator consisting of a mass $m$, a
linear spring of stiffness $k$, and a magnetorheological (MR) damper whose
damping coefficient is adjusted by an external magnetic field; see
Figure~\ref{fig:schematic}. The system is sampled at the instants
$t_k=kT$, and $x_k$ denotes the normalized displacement of the mass at the
$k$-th sampling instant. In the overdamped (quasi-static) regime the
one-step map of the normalized displacement is multiplicative,
\begin{equation}\label{eq:damper_map}
x_{k+1}=\lambda\bigl(|x_k|\bigr)\,x_k ,
\end{equation}
where $\lambda(|x_k|)\in(0,1)$ is the effective damping ratio over one
sampling period; the value of $\lambda$ is controlled through the magnetic
field. Indeed, for an overdamped oscillator the normalized displacement
satisfies $\tau(x)\dot x+x=0$ with a state-dependent relaxation time
$\tau>0$, so that after one sampling period
$x_{k+1}=e^{-T/\tau(x_k)}x_k$, which is \eqref{eq:damper_map} with
$\lambda(r)=e^{-T/\tau(r)}\in(0,1)$. Suppose that the controller modulates the field harmonically with
the measured displacement, so that
\begin{equation}\label{eq:damper_law}
\lambda(r)=\frac12+\frac13\sin(\kappa r),\qquad r\ge0,
\end{equation}
for some $\kappa>0$; after the normalization $\kappa=1$ the law is
\eqref{eq:harmonic}. The damping ratio satisfies
\begin{equation}\label{eq:damper_range}
\frac16\le\lambda(r)\le\frac56,\qquad r\ge0,
\end{equation}
so the damper never injects energy and the sampled map is never expansive
in the radial direction.

\begin{proposition}[Global asymptotic stability]\label{prop:damper}
For every $\kappa>0$ and every initial displacement $x_0\in\R$, the
solution of \eqref{eq:damper_map}--\eqref{eq:damper_law} satisfies
\begin{equation}\label{eq:damper_rate}
|x_k|\le\Bigl(\frac56\Bigr)^{k}|x_0|,\qquad k\ge0 .
\end{equation}
Consequently $x_k\to0$ for every $x_0$, and $0$ is the unique
equilibrium of the sampled system.
\end{proposition}

\begin{proof}
Put $\mu(r):=\frac12+\frac13\sin(\kappa r)$ for $r\ge0$. Then $\mu$ is
continuous and $0\le\mu(r)\le5/6$ for all $r\ge0$; moreover
$\mu(\norm{x})=\lambda(\norm{x})$ with $\lambda$ as in
\eqref{eq:damper_law}. Hence \eqref{eq:damper_map} is the radial damping
map \eqref{eq:T} associated with the law $\mu$, and
Theorem~\ref{thm:family} and Corollary~\ref{cor:radial_convergence}
apply with $\rho=5/6$ and yield \eqref{eq:damper_rate}, the convergence
$x_k\to0$, and the uniqueness of the equilibrium.
\end{proof}

\begin{figure}[htbp]
\centering
\includegraphics[width=0.88\textwidth]{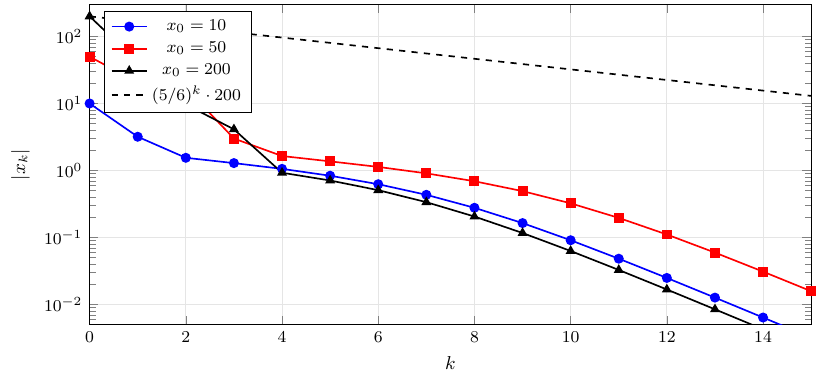}
\caption{Decay of the displacement under the harmonic law
\eqref{eq:harmonic}, for the initial displacements $x_0=10$, $50$ and
$200$, on a logarithmic scale; the dashed line is the worst-case bound
\eqref{eq:damper_rate}. All trajectories decay geometrically, although the
sampled map is not uniformly contracting.}
\label{fig:decay}
\end{figure}

Proposition~\ref{prop:damper} shows that the vibration is suppressed from
every initial displacement, with the geometric rate $5/6$; the bound is
sharp, because at the states with $\kappa r=\pi/2+2k\pi$ (where
$\sin(\kappa r)=1$) the one-step factor equals $5/6$ exactly; in the
normalized variables $\kappa=1$ these are the states
$r=\pi/2+2k\pi$. What is not true, and this is the
content of Section~\ref{sec:strict}, is that the sampled map is uniformly
contracting: by Theorem~\ref{thm:dispersion} and Example~\ref{ex:sin},
the map is not a Kirk asymptotic contraction. Physically, in the
normalized model, two states at distance $1$ but at large amplitude have
different damping ratios, and
their images can be separated by about $0.15\,r$ in one sampling period;
Table~\ref{tab:disp} and Figure~\ref{fig:disp} quantify this dispersion,
and Figure~\ref{fig:decay} shows the decay of $|x_k|$ for several initial
displacements.
The classical Kirk theorem, whose modulus is a function of the separation
alone, cannot be applied; the relaxed theorem of \cite{Shi2026} applies,
because $U_T$ records the size of the states and hence the scale at which
the damping acts.

The practical conclusion is the following. A state-dependent, periodic
modulation of the damping ratio is admissible as far as global asymptotic
stability is concerned: Proposition~\ref{prop:damper} guarantees it.
However, such a modulation rules out any uniform one-step contraction
certificate, and hence any argument based on Kirk's condition (for
instance, contraction-based observer or synchronization designs). If a
uniform contraction certificate is required, the modulation has to be
slowly varying, so that the map becomes a uniform contraction.

\subsection{Digital maturity gaps in small and medium-sized enterprises}
\label{subsec:gap}

In contrast with the physical scenario of the previous subsection, this
subsection concerns the management side of digital transformation: the
firm adjusts its investment to the gap between its digital maturity and
a target level. The adjustment is again multiplicative and state
dependent, so the two scenarios share the same mathematical structure;
what differs is the state, the adjustment rate and the management
interpretation.

Let $y_k\in[0,1]$ be the digital maturity of the firm in period $k$. In
practice the maturity index is a weighted composite of several
dimensions, for instance digital strategy and organization,
infrastructure and cloud services, data governance, digitalized business
processes, digital talent, and online channels and marketing; each
dimension is normalized to $[0,1]$ and the weights are nonnegative and
sum to one. Take the target level $1$ (or the industry benchmark) and
define the gap
\begin{equation}\label{eq:gap}
g_k:=1-y_k\ge0 .
\end{equation}
$g_k=0$ means that the target has been reached, while $g_k=1$ means that
the firm is not digitalized at all. The data may come from SME
digitalization assessments, platform questionnaires and annual reports,
usually at annual frequency.

Each period the firm invests a fraction of the gap:
\begin{equation}\label{eq:gap_map}
g_{k+1}=\bigl(1-\alpha(g_k)\bigr)g_k,
\end{equation}
where $\alpha(g_k)\in(0,1)$ is the fraction closed. The multiplicative
form means that resources are allocated in proportion to the remaining
gap, a simplification of the linear investment rules common in
management practice.

\begin{figure}[htbp]
\centering
\includegraphics[width=0.88\textwidth]{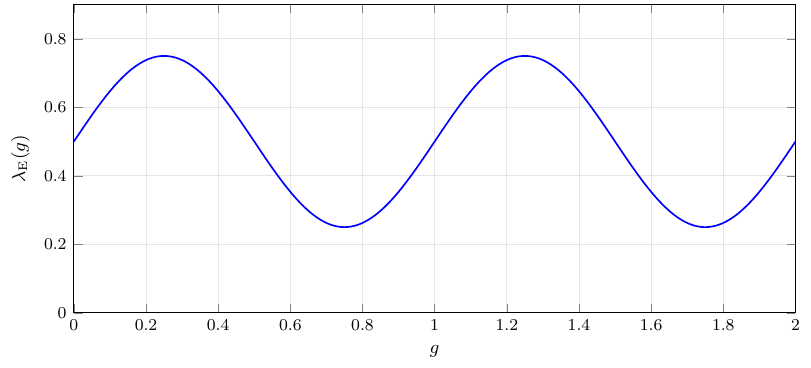}
\caption{The stage-gated adjustment law \eqref{eq:gap_law}: as the gap
varies within a stage, the adjustment factor oscillates periodically
between $1/4$ and $3/4$.}
\label{fig:gap_law}
\end{figure}

Stage-gated management makes the adjustment fraction depend on the
stage in which the gap currently lies. Measure the gap in digitalization
stages (one unit corresponds to one stage); within each stage the firm
goes through the cycle ``assessment--investment--consolidation'', and at
the stage boundaries (a new system goes live, new compliance or data
requirements appear) the adjustment intensity changes. As the simplest
smooth description take the first Fourier harmonic:
\begin{equation}\label{eq:gap_law}
\lambda(g):=1-\alpha(g)=\frac12+\frac14\sin(2\pi g),\qquad g\ge0 .
\end{equation}
Hence
\begin{equation}\label{eq:gap_range}
\frac14\le\lambda(g)\le\frac34,\qquad g\ge0,
\end{equation}
that is, at least $1/4$ and at most $3/4$ of the gap is closed each
period; the adjustment law is plotted in Figure~\ref{fig:gap_law}.

The form \eqref{eq:gap_law} can also be obtained from a simple
optimization problem. Suppose the firm chooses the fraction $u$ to
close, with benefit $b\,u\,g$ (efficiency and market gains from a
smaller gap) and cost $(\gamma/2)u^2$ (adjustment costs); without stage
effects the optimal choice is
\begin{equation}
u^*=\frac{b}{\gamma}\,g
\end{equation}
(capped at $1$). The fixed assessment costs and boundary effects of the
stage-gate process add a periodic component to the optimal fraction,
which after simplification gives the form \eqref{eq:gap_law}. We
emphasize that \eqref{eq:gap_law} is a stylized description of
stage-gated, state-dependent adjustment rather than an estimate for a
particular industry; the results below use only that $\lambda$ has a
positive lower bound, an upper bound strictly below one, and is
nonconstant and periodic.

\begin{proposition}[Convergence of the gap to zero]\label{prop:gap}
For every $g_0\ge0$ the solution of \eqref{eq:gap_map}--\eqref{eq:gap_law}
satisfies
\begin{equation}\label{eq:gap_rate}
0\le g_k\le\Bigl(\frac34\Bigr)^{k}g_0,\qquad k\ge0 .
\end{equation}
In particular $g_k\to0$ for every $g_0\ge0$, and $0$ is the unique
equilibrium in $[0,\infty)$.
\end{proposition}

\begin{proof}
Consider $T$ on $X=\R$. The law $\lambda$ is continuous with
\begin{equation}
1/4\le\lambda(g)\le3/4\qquad(g\ge0),
\end{equation}
so \eqref{eq:gap_map} is the radial damping map \eqref{eq:T} associated
with $\lambda$; moreover $T$ maps $[0,\infty)$ into itself because
$\lambda\ge0$. Theorem~\ref{thm:family} and
Corollary~\ref{cor:radial_convergence} apply with $\rho=3/4$ and give
\eqref{eq:gap_rate}, the convergence to $0$, and the uniqueness of the
equilibrium; the restriction to the invariant set $[0,\infty)$ preserves
these properties.
\end{proof}

Estimate \eqref{eq:gap_rate} yields the number of periods needed to
reach a given accuracy. For $\varepsilon>0$ one has $g_k\le\varepsilon$
as soon as
\begin{equation}
k\ge\frac{\ln(g_0/\varepsilon)}{\ln(4/3)}
\approx3.48\,\ln\frac{g_0}{\varepsilon}.
\end{equation}
For instance, halving the gap takes about $2.4$ periods, reducing it to
$10\%$ of its initial value takes about $8.0$ periods, and reducing it
to $1\%$ takes about $16.0$ periods. The actual trajectories are usually
faster, because $\lambda(0)=1/2$ and hence the adjustment rate
approaches $1/2$ when the gap is small. Table~\ref{tab:gap} reports
numerical trajectories for several initial gaps.

\begin{table}[htbp]
\centering
\begin{tabular}{cccccc}
\hline
$k$ & $g_0=0.2$ & $g_0=0.5$ & $g_0=1$ & $g_0=2$ &
$(3/4)^k\cdot2$ \\
\hline
$0$  & $0.2000$ & $0.5000$ & $1.0000$ & $2.0000$ & $2.0000$ \\
$2$  & $0.1033$ & $0.1875$ & $0.2500$ & $0.5000$ & $1.1250$ \\
$4$  & $0.0405$ & $0.0945$ & $0.1371$ & $0.1875$ & $0.6328$ \\
$6$  & $0.0122$ & $0.0359$ & $0.0605$ & $0.0945$ & $0.3559$ \\
$8$  & $0.0032$ & $0.0106$ & $0.0199$ & $0.0359$ & $0.2002$ \\
$10$ & $0.0008$ & $0.0028$ & $0.0055$ & $0.0106$ & $0.1126$ \\
\hline
\end{tabular}
\caption{Numerical trajectories for several initial gaps (iterating
\eqref{eq:gap_map}--\eqref{eq:gap_law}) and the worst-case bound; the
last column is the bound \eqref{eq:gap_rate} for $g_0=2$.}
\label{tab:gap}
\end{table}

The table also lists the worst-case bound $(3/4)^k\cdot2$ for $g_0=2$:
the actual trajectories always lie below the bound, and trajectories
corresponding to different initial gaps coincide after a few periods,
so the convergence path is relatively insensitive to the initial gap.

First, global convergence is guaranteed. As long as the fraction closed
per period has a positive lower bound (equivalently, $\lambda\le3/4$),
the gap converges geometrically to zero at rate $3/4$ for every initial
gap, that is, digital maturity converges to the target level. This
provides a theoretical basis for stage-gated, state-dependent adjustment
rules: a firm need not design a time-invariant linear rule in order to
obtain a stability guarantee.

Second, state-dependent rules produce transient divergence. Two firms
whose initial gaps differ by half a stage ($t=1/2$) can, within one
period, have their gap difference amplified to about $0.5$ times the
gap; hence single-period benchmarking data may show divergence across
firms even though all gaps converge in the long run. The management
implication is that benchmarking and policy evaluation should use
multi-period averages or comparisons within the same stage, rather than
single-period cross-sections.

Third, if a uniform convergence rate is required, a slowly varying
adjustment law should be used, so that the map becomes a uniform
contraction. Forecasting, cross-firm coordination and policy evaluation
often call for a rate independent of the state. There is a trade-off:
the stage-gated law has average adjustment rate $1/2$ and is faster on
average, but its worst case is $3/4$; the slowly varying law is slower on
average but guarantees a uniform rate.

Fourth, policy design can reduce short-run divergence by smoothing the
adjustment rate. Staged subsidies, funding calls and assessment deadlines
create periodic fluctuations of the adjustment rate; by
Proposition~\ref{prop:gap} such fluctuations do not destroy convergence,
but they affect the convergence path and short-run divergence.
Continuous funding, rolling applications and multi-year assessments can
reduce the amplitude of the adjustment rate and hence the transient
divergence.

Fifth, the management rule of thumb has two parts: ensure that the
fraction closed per period has a positive lower bound (avoid
$\alpha=0$ or negative values, i.e., no investment or digital
regression); and, if a uniform convergence certificate is required, use
a slowly varying adjustment law, whereas otherwise the stage-gated law
already guarantees global convergence.

Calibration requires firm-level panel data. The maturity index can be
composed from SME digitalization assessments, platform questionnaires or
annual reports; the fraction $\alpha(g)$ can be estimated by regressing
$y_{k+1}-y_k$ on $y_k$ nonparametrically or with Fourier terms, and by
testing whether the stage (periodic) component is significant. One
should account for the endogeneity of investment and the gap, measurement
error and sample selection; the model of this paper is illustrative
rather than an empirical estimate.

The model has three main limitations: it does not consider stochastic
shocks (demand fluctuations, technological change), firm heterogeneity
(different firms have different stage partitions and adjustment rates),
or entry and exit. Corresponding extensions include the convergence
speed and the degree of divergence under random perturbations, the
aggregation of heterogeneous firms, and the integration with policy
evaluation frameworks.

\section{Conclusion}\label{sec:discussion}

We have shown that the relaxed class of asymptotic contractions
introduced in \cite{Shi2026} strictly contains Kirk's class. The
separating examples are radial damping maps
$T(x)=\lambda(\norm{x})x$ with $0\le\lambda\le\rho<1$ whose one-step
dispersion is unbounded; explicit instances are the periodic harmonic law,
the slowly modulated law, and the two-dimensional harmonic law. These maps
are defined by a one-line formula, they are physically meaningful as
state-dependent damping laws, and they occur in every normed space. The
obstruction that destroys Kirk's condition is the unbounded one-step
dispersion \eqref{eq:dispersion_inf} at a fixed distance, and the
condition survives because the upper quasi-metric $U_T$ records the size
of the states, hence the scale of the contraction. On the other hand, in a
complete space, if the image $T(X)$ is relatively compact, the relaxed
condition collapses to Kirk's condition
(Theorem~\ref{thm:compact}), so the strictness is a non-proper,
``infinite-dimensional'' phenomenon: on compact or finite metric spaces
the two classes coincide.

Several questions remain open.
\begin{enumerate}
\item Is relative compactness of $T(X)$ necessary for the coincidence of
the two classes? Theorem~\ref{thm:compact} gives sufficiency, and
Corollary~\ref{cor:compact} shows that the question is nontrivial only in
spaces in which bounded sets need not be relatively compact.
\item Can strictness be witnessed by a map for which all the sup-moduli
\begin{equation}\label{eq:Mn}
M_n(t):=\sup\bigl\{d(T^nx,T^ny): d(x,y)\le t\bigr\},
\qquad t>0,\ n\ge0,
\end{equation}
are finite? In Example~\ref{ex:sin} the strictness is caused by
$M_1(1)=+\infty$; a witness with finite moduli would show that the
strictness is more delicate than the mere unboundedness of the dispersion
at a single scale.
\item Does the relaxed condition admit a Meir--Keeler (Suzuki type)
version, in which the Boyd--Wong condition on $\psi$ is replaced by a
Meir--Keeler condition, and is the resulting class strictly larger than
Suzuki's class of asymptotic contractions of Meir--Keeler type?
\item Can the convergence in Theorem~\ref{thm:quoted} be quantified in
terms of the moduli and of $d(x,z)$? For the radial damping maps the rate
is geometric, see Corollary~\ref{cor:radial_convergence}; for general
relaxed asymptotic contractions no such estimate is available.
\item Is the class of relaxed asymptotic contractions stable under small
perturbations of the map? In the radial family, a perturbation of the
damping law which remains in $[0,\rho]$ preserves the constants
$2\rho^{n+1}(1-\rho)^{-1}$, but a general stability theory is missing.
\end{enumerate}

\section*{Use of Generative-AI tools declaration}

The author declares that Artificial Intelligence (AI) tools were used for
language editing and for the preparation of the \LaTeX{} file of the
manuscript; the mathematical results and their proofs were checked by the
author.

\section*{Conflict of interest}

The author declares there is no conflict of interest.

\section*{Acknowledgments}

The author thanks the Department of Mathematics and Statistics of Hubei
Engineering University for its support.

\end{document}